\documentclass[12pt]{article}

\usepackage{amssymb,amsmath}
\usepackage{theorem}
\usepackage{pifont}

\newcommand {\nF} {\mathcal{F}}

\newcommand {\art}[6]{{\sc #1:} {#2.} {\em #3} {\bf #4} {(#5),} {#6.}}

\newcommand {\book}[5]{{\sc #1:} {``#2."} {#3,} {#4} {(#5).}}

\newcommand {\toappearin}[3]{{\sc #1:} {#2,}
                    {\sf{to appear in}} {\em #3}.}

\newcommand \func {{}^\omega\omega}

\newcommand \la {{\langle}}
\newcommand \ra {{\rangle}}
\newcommand \subs {{\subseteq}}
\newcommand \lesseq {\preccurlyeq}
\newcommand \nlesseq {\not\preccurlyeq}
\newcommand \card {\mathfrak{l}}
\newcommand {\nl} {\\[0.8ex]}
\newcommand \cf {{\sl cf.}\;}
\newcommand \eg {{\sl e.g.}\;}

\newcommand {\forces}[1] {{\hspace{0.5ex}}{\rule{0.1ex}{1.5ex}}
                    {\hspace{0.2ex}}{\rule{0.1ex}{1.5ex}}
                    {\rule[0.75ex]{1.5ex}{0.1ex}}
                    {\hspace{0.2ex}}{}_{#1}{\hspace{0.4ex}}}

\newcommand \cP {{\mathcal{P}}}

\newcommand \R {\mathbb{R}}
\newcommand \MM {{\mathbb M}\mskip1mu{}}
\newcommand \Mi {{\mathbb M}\mskip1mu{}_{\omega_2}}

\newcommand \fp {{\mathfrak{p}}}
\newcommand \fb {{\mathfrak{b}}}
\newcommand \fd {{\mathfrak{d}}}
\newcommand \fc {{\mathfrak{c}}}

\newcommand \mmin {\operatorname{min}}

\newcommand \ZFC {{\text{\sf ZFC}}}
\newcommand \MA {{\text{\sf MA}}}
\newcommand \CH {{\text{\sf CH}}}

\def\phi{\varphi}
\def\epsilon{\varepsilon}

\theorembodyfont{\rmfamily}
\newtheorem {nummer}{ }[section]
\newtheorem {thm}[nummer]{{\sc{Theorem}}\rm\small}
\newtheorem {prop}[nummer]{{\sc{Proposition}}\rm\small}
\newtheorem {lm}[nummer]{{\sc{Lemma}}\rm\small}
\newtheorem {fct}[nummer]{{\sc{Fact}}\rm\small}
\newtheorem {cor}[nummer]{{\sc{Corollary}}\rm\small}
\newcommand \rmk {{\sc Remark:}\hspace*{3mm}}

\newcommand \problem  {{\sc Problem:}\hspace*{3mm}}
\newcommand \proof {\noindent{\sc Proof:}\hspace*{3mm}}

\def\eop{{\unskip\nobreak\hfil\penalty50\hskip8mm\hbox{}
  \nobreak\hfil
  {$\boldsymbol{\dashv}$}\parfillskip=0mm \par\smallskip}}

\begin{document}

\begin{center}
      {\large{\bf{ON A THEOREM OF BANACH AND KURATOWSKI}}}\\[2ex]
      {\large{\bf{AND $\boldsymbol{K}$-LUSIN SETS}}}
\end{center}
\smallskip

\begin{center}
{\sc{Tomek Bartoszy\'{n}ski}}\\[1.5ex] {\small{\sl Department of
Mathematics, Boise State University\\ Boise, ID 83725
(U.S.A.)}}\\[0.5ex] {\small{\sl Email: tomek@diamond.boisestate.edu}}
\end{center}

\begin{center}
{\sc{Lorenz Halbeisen}}\\[1.5ex] {\small{\sl Department of Pure
Mathematics, Queen's University Belfast\\ Belfast BT7 1NN (Northern
Ireland)}}\\[0.5ex] {\small{\sl Email: halbeis@qub.ac.uk}}
\end{center}
\medskip

\begin{abstract}\noindent
In a paper of 1929, Banach and Kuratowski proved---assuming the
continu\-um hypothesis---a combinatorial theorem which implies that
there is no non-vanishing $\sigma$-additive finite measure $\mu$ on
$\R$ which is defined for every set of reals. It will be shown that
the combinatorial theorem is equivalent to the existence of a
$K$-Lusin set of size $\mathfrak{2}^{\aleph_0}$ and that the
existence of such sets is independent of $\ZFC+\neg\CH$.
\end{abstract}
\renewcommand{\thefootnote}{}
\footnotetext{{}\hfill\\[-1.5ex] {\it 2000 Mathematics Subject
Classification:} {\bf 03E35} 03E17 03E05
First author  partially supported by
NSF grant DMS 9971282 and Alexander von Humboldt Foundation}
\setcounter{section}{-1}
\section{Introduction}

In \cite{BanachKuratowski}, Stefan Banach and Kazimierz Kuratowski
investigated the following problem in measure theory:

\noindent \problem {\it Does there exist a non-vanishing finite
measure $\mu$ on $[0,1]$ defined for every $X\subs [0,1]$, which is
$\sigma$-additive and such that for each $x\in [0,1]$,
$\mu(\{x\})=0\,$?}

They showed that such a measure does not exist if one assumes the
continuum hypothesis, denoted by $\CH$. More precisely, assuming
$\CH$, they proved a combinatorial theorem
(\cite[Th\'eor\`eme\,II]{BanachKuratowski}) and showed that this
theorem implies the non-existence of such a measure. The
combinatorial result is as follows:

\noindent {\sc Banach-Kuratowski Theorem\,} Under the assumption of
$\CH$, there is an infinite matrix $A^i_k\subs [0,1]$ (where
$i,k\in\omega$) such that:\vspace{-1.6ex}
\begin{itemize}
\item[{(i)}] For each $i\in\omega$,
$[0,1]=\bigcup_{k\in\omega}A^i_k$.
\item[{(ii)}] For each $i\in\omega$, if $k\neq k'$ then $A^i_k\cap
A^i_{k'}=\emptyset$.
\item[{(iii)}] For every sequence $k_0,k_1,\ldots,k_i,\ldots$ of
$\omega$, the set $\bigcap_{i\in\omega}(A^i_{0}\cup A^i_{1}\cup
\ldots\cup A^i_{k_i})$ is at most countable.
\end{itemize}

In the sequel, we call an infinite matrix $A^i_k\subs [0,1]$ (where
$i,k\in\omega$) for which (i), (ii) and (iii) hold, a {\bf
BK-Matrix}.

Wac{\l}aw Sierpi\'nski proved---assuming $\CH$---in
\cite{Sierpinski29} and \cite{Sierpinski32} two theorems involving
sequences of functions on $[0,1]$, and showed in \cite{Sierpinski29}
and \cite{Sierpinski33} that these two theorems are equivalent to the
Banach-Kuratowski Theorem, or equivalently, to the existence of a
BK-Matrix.

\noindent \rmk Concerning the problem in measure theory mentioned
above, we like to recall the well-known theorem of Stanis{\l}aw Ulam
(\cf \cite{Ulam} or \cite[Theorem\,5.6]{Oxtoby}), who showed that
each $\sigma$-additive finite measure $\mu$ on $\omega_1$, defined
for every set $X\subs\omega_1$ with $\mu(\{x\})=0$ for each
$x\in\omega_1$, vanishes identically. This result implies that if
$\CH$ holds, then there is no non-vanishing $\sigma$-additive finite
measure on $[0,1]$.

In the sequel we show that even if $\CH$ fails, the existence of a
BK-Matrix---which will be shown to be equivalent to the existence of
a $K$-Lusin set of size $\mathfrak{2}^{\aleph_0}$---may still be
true.

Our set-theoretical terminology (including forcing) is standard and
may be found in textbooks like \cite{BartoszynskiJudah},
\cite{Jechbook} and \cite{Kunen}.

\section{The Banach-Kuratowski Theorem revisited}

Before we give a slightly modified version of the Banach-Kuratowski
proof of their theorem, we introduce some notation.

For two functions $f,g\in\func$ let $f\lesseq g$ if and only if for
each $n\in\omega$, $f(n)\le g(n)$.

For $\nF\subs\func$, let $\lambda(\nF)$ denote the least cardinality
$\lambda$, such that for each $g\in\func$, the cardinality of
$\{f\in\nF: f\lesseq g\}$ is strictly less than $\lambda$. If
$\nF\subs\func$ is a family of size $\fc$, where $\fc$ is the
cardinality of the continuum, then we obviously have
$\aleph_1\le\lambda(\nF)\le\fc^+$. This leads to the following
definition: $$\card :=\mmin\{\lambda(\nF): \nF\subs\func\wedge
|\nF|=\fc\}\,.$$ If one assumes $\CH$, then one can easily construct
a family $\nF\subs\func$ of cardinality $\fc$ such that
$\lambda(\nF)=\aleph_1$, and therefore, $\CH$ implies that $\card
=\aleph_1$.

The crucial point in the Banach-Kuratowski proof of their theorem is
\cite[Th\'eor\`eme\,II']{BanachKuratowski}. In our notation, it reads
as follows:

\begin{prop}\label{prop:equiv}
The existence of a BK-Matrix is equivalent to $\card =\aleph_1$.
\end{prop}

For the sake of completeness and for the reader's convenience, we
give the Banach-Kuratowski proof of Proposition\,\ref{prop:equiv}.

\proof $(\Leftarrow)\,$ Let $\nF\subs\func$ be a family of
cardinality $\fc$ with $\lambda(\nF)=\aleph_1$. In particular, for
each $g\in\func$, the set $\{f\in\nF: f\lesseq g\}$ is at most
countable. Let $f_\alpha$ ($\alpha<\fc$) be an enumeration of $\nF$.
Since the interval $[0,1]$ has cardinality $\fc$, there is a
one-to-one function $\Xi$ from $[0,1]$ onto $\nF$. For $x\in [0,1]$,
let $n^x_i:=\Xi(x)(i)$. Now, for $i,k\in\omega$, define the sets
$A^i_k\subs [0,1]$ as follows: $$x\in A^i_k\ \text{if and only if}\
k=n^x_i\,.$$ It is easy to see that these sets satisfy the conditions
(i) and (ii) of a BK-Matrix. For (iii), take any sequence
$k_0,k_1,\ldots,k_i,\ldots$ of $\omega$ and pick an arbitrary
$x\in\bigcap_{i\in\omega}(A^i_0\cup A^i_1\cup \ldots\cup A^i_{k_i})$.
By definition, for each $i\in\omega$, $x$ is in $A^i_0\cup A^i_1\cup
\ldots\cup A^i_{k_i}$. Hence, for each $i\in\omega$ we get $n^x_i\le
k_i$, which implies that for $g\in\func$ with $g(i):=k_i$ we have
$\Xi(x) \lesseq g$. Now, since $\lambda(\nF)=\aleph_1$,
$\Xi(x)\in\nF$ and $x$ was arbitrary, the set $\{x\in [0,1]:
\Xi(x)\lesseq g\}=\bigcap_{i\in\omega}(A^i_0\cup A^i_1\cup \ldots\cup
A^i_{k_i})$ is at most countable.\\[0.7ex] $(\Rightarrow)\,$ Let
$A^i_k\subs [0,1]$ (where $i,k\in\omega$) be a BK-Matrix and let
$\nF\subs\func$ be the family of all functions $f\in\func$ such that
$\bigcap_{i\in\omega}A^i_{f(i)}$ is non-empty. Is is easy to see that
$\nF$ has cardinality $\fc$. Now, for any sequence
$k_0,k_1,\ldots,k_i,\ldots$ of $\omega$, the set
$\bigcap_{i\in\omega}(A^i_0\cup A^i_1\cup \ldots\cup A^i_{k_i})$ is
at most countable, which implies that for $g\in\func$ with
$g(i):=k_i$, the set $\{f\in\nF: f\lesseq g\}$ is at most countable.
Hence, $\lambda(\nF)=\aleph_1$. \eop

\section{$\boldsymbol{K}$-Lusin sets}\label{sec:Lusin}

In this section we show that $\card =\aleph_1$ is equivalent to the
existence of a $K$-Lusin set of size $\fc$.

We work in the Polish space $\func$.

\begin{fct}\label{fct:compact}
A closed set $K\subs\func$ is compact if and only if there is a
function $f\in\func$ such that $K\subs\{g\in\func: g\lesseq f\}$.
\end{fct}

(See \cite[Lemma\,1.2.3]{BartoszynskiJudah}.)

An uncountable set $X\subs\func$ is a {\bf Lusin set}, if for each
meager set $M\subs\func$, $X\cap M$ is countable.

An uncountable set $X\subs\func$ is a {\bf $\boldsymbol{K}$-Lusin
set}, if for each compact set $K\subs\func$, $X\cap K$ is countable.

\begin{lm}\label{lm:compact-meager} Every Lusin set is a $K$-Lusin set.
\end{lm}

\proof By Fact\,\ref{fct:compact}, every compact set $K\subs\func$ is
meager (even nowhere dense), and therefore, every Lusin set is a
$K$-Lusin set. \eop

\begin{lm}\label{lm:equiv} The following are equivalent:\\
\hspace*{2ex}{(a)} $\card =\aleph_1$.\\ \hspace*{2ex}{(b)} There is a
$K$-Lusin set of cardinality $\fc$.
\end{lm}

\proof This follows immediately from the definitions and
Fact\,\ref{fct:compact}. \eop

\noindent \rmk Concerning Lusin sets we like to mention that
Sierpi\'nski gave in \cite{Sierpinski37} a combinatorial result which
is equivalent to the existence of a Lusin set of cardinality $\fc$.

For $f,g\in\func$, define $f\lesseq^*g$ if $f(n)\le g(n)$ for all but
finitely many $n\in\omega$. The cardinal numbers $\fb$ and $\fd$ are
defined as follows: $$\fb:=\mmin\big{\{}|\nF|: \nF\subs\func\
\text{{\sl and\/}}\ \forall g\in\func\exists
f\in\nF\;(f\nlesseq^*g)\big{\}}$$ $$\fd:=\mmin\big{\{}|\nF|:
\nF\subs\func\ \text{{\sl and\/}}\ \forall g\in\func\exists
f\in\nF\;(g\lesseq^*f)\big{\}}$$

\begin{lm}\label{lm:bd} $\card =\aleph_1$ implies $\fb=\aleph_1$ and
$\fd=\fc$. Moreover, $K$-Lusin sets are exactly those (uncountable)
subsets of $\func$ whose all uncountable subsets are unbounded.
(Families like that are also called {\it strongly unbounded} and they
play an important role in preserving unbounded families in
iterations, see \eg \cite{BartoszynskiJudah} for details.)
\end{lm}

\proof Assume $\card =\aleph_1$, then, by Lemma\,\ref{lm:equiv},
there exists a $K$-Lusin set $X\subs\func$ of cardinality $\fc$. It
is easy to see that every uncountable subset of $X$ is unbounded, so,
$\fb=\aleph_1$. On the other hand, every function $g\in\func$
dominates only countably many elements of $X$. Hence, no family
$\nF\subs\func$ of cardinality strictly less than $\fc$ can dominate
all elements of $X$, and thus, $\fd=\fc$. \eop

\begin{prop}\label{prop:Lusin}
Adding $\kappa$ many Cohen reals produces a Lusin set of size
$\kappa$.
\end{prop}

(See \cite[Lemma\,8.2.6]{BartoszynskiJudah}.)

\begin{thm}\label{thm:main} The existence of a $K$-Lusin set
of cardinality $\fc$ is independent of $\ZFC+\neg\CH$.
\end{thm}

\proof By Proposition\,\ref{prop:Lusin} and
Lemma\,\ref{lm:compact-meager} it is consistent with $\ZFC$ that
there is a $K$-Lusin set of cardinality $\fc$.\\ On the other hand,
it is consistent with $\ZFC$ that $\fb>\aleph_1$ or that $\fd<\fc$
(\cf \cite{BartoszynskiJudah}). Therefore, by Lemma\,\ref{lm:bd}, it
is consistent with $\ZFC$ that there are no $K$-Lusin sets of
cardinality $\fc$. \eop

By Lemma\,\ref{lm:equiv} and Proposition\,\ref{prop:equiv}, as a
immediate consequence of Theorem\,\ref{thm:main} we get the
following.

\begin{cor} The existence of a BK-Matrix is independent of
$\ZFC+\neg\CH$.
\end{cor}

\section{Odds and Ends}

An uncountable set $X\subs [0,1]$ is a {\bf Sierpi\'{n}ski set}, if
for each measure zero set $N\subs [0,1]$, $X\cap N$ is countable.

\begin{prop}\label{prop:rothberger}
The following are equivalent:\nl \hspace*{2ex}{(a)} $\CH$.\\
\hspace*{2ex}{(b)} There is a Lusin set of cardinality $\fc$ and an
uncountable Sierpi\'{n}ski set.\\ \hspace*{2ex}{(c)} There is a
Sierpi\'{n}ski set of cardinality $\fc$ and an uncountable Lusin set.
\end{prop}

(See \cite[p.\,217]{Rothberger}.)

\begin{prop} It is consistent with $\ZFC$ that there exists a
$K$-Lusin set of cardinality $\fc$, but there are neither Lusin nor
Sierpi\'{n}ski sets.
\end{prop}

\proof Let $\Mi$ denote the $\omega_2$-iteration of Miller
forcing---also called ``rational perfect set forcing''---with
countable support. Let us start with a model $V$ in which $\CH$ holds
and let $G_{\omega_2}=\la m_\iota:\iota<\omega_2\ra$ be the
corresponding generic sequence of Miller reals. Then, in
$V[G_{\omega_2}]$, $G_{\omega_2}$ is a $K$-Lusin set of cardinality
$\fc=\aleph_2$. For this, we have to show the following property:
$$\text{\sl For all\/}\ f\in\func\cap V[G_{\omega_2}],\ \text{\sl the
set\/}\ \{\iota: m_\iota\lesseq f\}\ \text{\sl is countable.}$$
Suppose not and let $f\in\func\cap V[G_{\omega_2}]$ be a witness.
Further, let $p$ be an $\Mi$-condition such that $$p\forces{\Mi}
``\text{\sl for some $n_0\in\omega$, the set\/}\ \big{\{}\iota:
\forall k\ge n_0 \big(m_\iota(k) < \dot{f}(k)\big)\big{\}}\ \text{\sl
is uncountable''.}$$ We can assume that these dominated reals
$m_\iota$ are among $\{m_\alpha: \alpha < \beta <\omega_2\}$ and that
$\beta$ is minimal. This way, $f$ is added after step $\beta$ of the
iteration. Let $a^*:=\operatorname{cl}(\dot{f})$ be the (countable)
set of ordinals such that if we know $\{m_\iota: \iota \in a\}$, then
we can compute $\dot{f}$. (Notice that $a^*$ is much more than just
the support of $\dot{f}$, since it contains also supports of all
conditions that are involved in conditions involved in $\dot{f}$, and
so on.) Let $N$ be a countable model such that $p,\dot{f}\in N$,
$a^*\subs N$ and let $\MM_{a^*}$ be the iteration of Miller forcing,
where we put the empty forcing at stages $\alpha \notin a^*$
(essentially, $\MM_{a^*}$ is the same as
$\MM_{\operatorname{o.t.}(a^*)}$).\\
 The crucial lemma---which is done in \cite[Lemma\,3.1]{ShelahSpinas}
for Mathias forcing, but also works for Miller forcing---is the
following: If $N \models p\in \MM_{a^*}$, then there exists a
$q\in\Mi$ which is stronger than $p$ such that
$\operatorname{cl}(q)={a^*}$ and $q$ is $(N,\MM_{a^*})$-generic over
$N$. In particular, if $\{m_\iota
:\iota<\omega_2\}$ is a generic sequence of Miller reals consistent
with $q$, then $\{m_\iota : \iota\in {a^*}\}$ is $\MM_{a^*}$-generic
over $N$ (consistent with $p$).\\
 So, fix such a $q$. Now we claim that for $\gamma\in\beta\setminus
N$, $q$ forces that $\dot{f}(k)> m_\gamma(k)$ for some $k\ge n_0$:
Take any $\gamma\in\beta\setminus N$ and let $q^*$ be a condition
stronger than $q$. Let $q^*_1=q_1 | \beta$, and let $q^*_2= q^* | a$.
Without loss of generality, we may assume that $q^*_2 =q$. Now, first
we strengthen $q^*_1$ to determine the length of stem of
$q^*_1(\gamma)$, and make it equal to some $k>n_0$. Next we shrink
$q^*_2$ to determine the first $k$ digits of $\dot{f}$. Finally, we
shrink $q^*_1(\gamma)$ such that $q^*_1(\gamma)(k)>\dot{f}(k)$. Why
we can do this? Although $f$ is added after $m_\gamma$, from the
point of view of model $N$, it was added before. So, working below
condition $q^*_2$ (in $\MM_{a^*}$) we can compute as many digits of
$\dot{f}$ as we want without making any commitments on $m_\gamma$,
and vice versa. Even though the computation is in $N$, it is
absolute. This completes the first part of the proof.\\
 On the other hand, it is known (\cf \cite{JudahShelah}) that in
$V[G_{\omega_2}]$, there are neither Lusin nor Sierpi\'nski sets of
any uncountable size, which completes the proof. \eop

\begin{prop} It is consistent with $\ZFC$ that $\fb=\aleph_1$ and
$\fd=\fc$, but there is no $K$-Lusin set of cardinality $\fc$.
\end{prop}

\proof Take a model $M$ in which we have $\fc=\aleph_2$ and in which
Martin's Axiom $\MA$ holds. Let $G=\la c_\beta:\beta<\omega_1\ra $ be
a generic sequence of Cohen reals of length $\omega_1$. In the
resulting model $M[G]$ we have $\fb=\aleph_1$ (since the set of Cohen
reals forms an unbounded family) and $\fd=\aleph_2$. On the other
hand, there is no $K$-Lusin set of cardinality $\fc$ in $M[G]$. Why?
Suppose $X\subs\func$ has cardinality $\aleph_2$. Take a countable
ordinal $\alpha$ and a subset $X'\subs X$ of cardinality $\aleph_2$
such that $X'\subs M[G_\alpha]$, where $G_\alpha:=\la
c_\beta:\beta\le\alpha\ra$. Now, $M[G_\alpha]=M[c]$ (for some Cohen
real $c$) and $M[c]\models \MA(\sigma\text{\sl -centered\/})$, (\cf
\cite{Roitman} or \cite[Theorem\,3.3.8]{BartoszynskiJudah}). In
particular, since $\MA(\sigma\text{\sl -centered\/})$ implies
$\fp=\fc$ and $\fp\le\fb$, we have $M[c]\models\fb=\aleph_2$. Thus,
there is a function which bounds uncountably many elements of $X'$.
Hence, by Lemma\,\ref{lm:bd}, $X$ cannot be a $K$-Lusin set. \eop

Let $Q$ be a countable dense subset of the interval $[0,1]$. Then
$X\subs [0,1]$ is {\bf concentrated on $Q$}, if every open set of
$[0,1]$ containing $Q$, contains all but countably many elements of
$X$.

\begin{prop} The following are equivalent:\nl
\hspace*{2ex}{(a)} There exists a $K$-Lusin set of cardinality
$\fc$.\\ \hspace*{2ex}{(b)} There exists a concentrated set of
cardinality $\fc$.
\end{prop}

\proof (b)$\rightarrow$(a)\quad Let $Q$ be a countable dense set in
$[0,1]$ and let $\phi:[0,1]\setminus Q\to\func$ be a homeomorphism.
If $U\subs\func$ is compact, then $\phi^{-1}[U]$ is compact, so
closed in $[0,1]$, and $[0,1]\setminus\phi^{-1}[U]$ is an open set
containing $Q$. Hence, the image under $\phi$ of an uncountable set
$X\subs [0,1]$ concentrated on $Q$ is a $K$-Lusin set of the same
cardinality as $X$.\\ (a)$\to$(b)\quad The preimage under $\phi$ of a
$K$-Lusin set of cardinality $\fc$ is a set concentrated on $Q$ of
the same cardinality. \eop

\noindent \rmk A Lusin set is concentrated on every countable dense
set, and concentrated sets have always strong measure zero. However,
the existence of a strong measure zero set of size $\fc$ does not
imply the existence of a concentrated sets of size $\fc$. In fact,
the existence of a strong measure zero set of size $\fc$ is
consistent with $\fd =\aleph_1$ (see \cite{BartoszynskiShelah}).

\end{document}